\documentclass[11.05pt,a4paper]{amsart}
\usepackage{amssymb,amsmath}
\usepackage{color}
\usepackage{latexsym}
\usepackage{amsthm,amsfonts,amssymb,mathrsfs}
\usepackage{rotating}
\usepackage[leqno]{amsmath}
\usepackage{xspace}
\usepackage[all]{xy}
\usepackage{longtable}

\headsep=1cm \numberwithin{equation}{section}
\newtheorem{theorem}{Theorem}[section]
\newtheorem{lemma}[theorem]{Lemma}
\newtheorem{proposition}[theorem]{Proposition}
\newtheorem{corollary}[theorem]{Corollary}

\theoremstyle{definition}
\newtheorem{definition}[theorem]{Definition}
\theoremstyle{remark}
\newtheorem{remark}[theorem]{Remark}
\newtheorem{example}[theorem]{Example}

\newcommand{\Mod}{\operatorname{Mod}}
\newcommand{\im}{\operatorname{im}}
\newcommand{\Ker}{\operatorname{Ker}}

\newcommand{\id}{\operatorname{id}}

\newcommand{\Gpd}{\operatorname{Gpd}}

\newcommand{\Ext}{\operatorname{Ext}}

\newcommand{\Hom}{\operatorname{Hom}}

\newcommand{\Coker}{\operatorname{Coker}}

\newcommand{\BZ}{\Bbb Z}

\newcommand{\lo}{\longrightarrow}
\newcommand{\fm}{\frak{m}}

\newenvironment{prf}[1][Proof]{\begin{proof}[\bf #1]}{\end{proof}}
\begin{document}

\author[F. Zareh-Khoshchehreh, M. Asgharzadeh and K. Divaani-Aazar]{Fatemeh Zareh-Khoshchehreh, Mohsen Asgharzadeh
and  Kamran Divaani-Aazar}
\title[Gorenstein homology, relative pure homology and ...]
{Gorenstein homology, relative pure homology and virtually Gorenstein rings}

\address{F. Zareh-Khoshchehreh, Department of Mathematics, Alzahra University, Vanak, Post Code 19834, Tehran,
Iran.}
\email{fzarehkh@gmail.com}

\address{M. Asgharzadeh, School of Mathematics, Institute for Research in Fundamental Sciences (IPM), P.O. Box
19395-5746, Tehran, Iran.}
\email{asgharzadeh@ipm.ir}

\address{K. Divaani-Aazar, Department of Mathematics, Alzahra University, Vanak, Post Code 19834, Tehran, Iran- and-School of Mathematics, Institute for Research in Fundamental Sciences (IPM), P.O. Box 19395-5746, Tehran, Iran.}
\email{kdivaani@ipm.ir}


\keywords {Cotorsion theories; Gorenstein flat modules; Gorenstein injective modules; Gorenstein injective preenvelopes; Gorenstein projective modules; Gorenstein projective precovers;  pure flat modules; pure injective modules, pure projective modules; virtually Gorenstein rings.\\
The research of the second and third authors are supported by grants from
IPM (No. 91130407 and No. 91130212, respectively).}

\begin{abstract} We consider the following question: Is Gorenstein homology a $\mathscr{X}$-pure homology,
in the sense defined by Warfield, for a class $\mathscr{X}$ of modules? Let $\mathcal{GP}$ denote the class 
of Gorenstein projective modules. We prove that over a commutative Noetherian ring $R$ of finite Krull 
dimension, Gorenstein homology is a $\mathcal{GP}$-pure homology if and only if $R$ is virtually Gorenstein.
\end{abstract}

\maketitle

\section{Introduction}

Throughout this paper, $R$ will be a commutative ring with identity. Let $\mathcal{GP}$ and $\mathcal{GI}$
denote the classes of Gorenstein projective and Gorenstein injective $R$-modules, respectively. We refer the
reader to \cite{AB}, \cite{EJ1}, \cite{EJ2} and \cite{EJT} for all unexplained definitions in the sequel.
Warfield \cite{W} has introduced a notion of $\mathscr{X}$-purity for any class $\mathscr{X}$ of $R$-modules.
Recall from \cite{W} that an exact sequence $0\to A\to B\to C\to 0$ of $R$-modules is called $\mathscr{X}$-
{\em pure exact} if for any $U\in \mathscr{X}$ the induced $R$-homomorphism $\Hom_{R}(U,B)\to \Hom_{R}(U,C)$
is surjective. An $R$-module $M$ is called $\mathscr{X}$-{\em pure projective} (respectively $\mathscr{X}$-
{\em pure injective}; $\mathscr{X}$-{\em pure flat}) if the functor $\Hom_R(M,-)$ (respectively $\Hom_R(-,M)$;  $M\otimes_R-$) leaves any $\mathscr{X}$-pure exact sequence exact. For a survey of results on purity, we refer
the reader to \cite{Hu}, \cite{S} and \cite{W}.

Let $\mathscr{X}$ be a class of $R$-modules. We say a homology theory $\mathcal{T}$ is a $\mathscr{X}$-{\it 
pure homology} if an $R$-module $M$ is projective (respectively injective; flat) in $\mathcal{T}$ if and only
if it is $\mathscr{X}$-pure projective (respectively $\mathscr{X}$-pure injective; $\mathscr{X}$-pure flat).
In this paper, we investigate the question: Is Gorenstein homology a $\mathscr{X}$-pure homology for an
appropriate class $\mathscr{X}$ of $R$-modules? Our candidate for a such class $\mathscr{X}$ is $\mathcal{GP}$.
To treat this question, we focus on Noetherian rings of finite Krull dimension. So, from now to the end of
the introduction, assume that $R$ is Noetherian of finite Krull dimension. It is not hard to verify that $$\mathcal{GP}=\{M\in R-\Mod|M\  \text{is} \  \mathcal{GP}-\text{pure projective}\};$$ see Proposition \ref{26}
and Remark \ref {36} (iv) below.  On the other hand,  we show that if the classes of Gorenstein injectives and $\mathcal{GP}$-pure injectives are the same, then also the classes of Gorenstein flats and $\mathcal{GP}$-pure 
flats are the same; see Corollary \ref{32} below. Therefore, our question reduces to: Are the classes of 
Gorenstein injectives and $\mathcal{GP}$-pure injectives the same? In Lemma \ref {38}, we show that $^{\bot}\mathcal{GI}\subseteq \mathcal{GP}^{\bot};$  see the beginning of the next section for the definitions of $^{\bot}\mathcal{GI}$ and $\mathcal{GP}^{\bot}.$ We call $R$ \emph{virtually Gorenstein} if $^{\bot}\mathcal{GI}=\mathcal{GP}^{\bot}.$ This generalizes the notion of virtually Gorenstein Artin algebras which was introduced by Beligiannis and Reiten in \cite{BR}. Our main results make such algebras relevant also in commutative ring theory. We prove that Gorenstein homology is a $\mathcal{GP}$-pure homology if and only if the functor $\Hom_R(-,\sim)$ is right balanced by $\mathcal{GP}\times\mathcal{GI}$ and if and only if $R$ is virtually Gorenstein; see Theorems \ref{41} and
\ref{43} below.

\section{Gorenstein And Relative Pure Projectivity}

Recall from the introduction that $\mathcal{GP}$ and $\mathcal{GI}$ denote the classes of Gorenstein projective
and Gorenstein injective $R$-modules, respectively. Our main result in this section is that when $\mathcal{GP}$
is precovering, then  $$\mathcal{GP}=\{M\in R-\Mod|M \  \  \text{is} \  \  \mathcal{GP}-\text{pure projective}\}.$$

For any class $\mathscr{X}$ of $R$-modules, let $\mathscr{X}^{\bot}$ (respectively $^{\bot}\mathscr{X}$) denote
the class of $R$-modules $M$ with the property that $\Ext^{1}_{R}(X,M)=0$ (respectively $\Ext^{1}_{R}(M,X)=0$)
for all $R$-modules $X\in \mathscr{X}.$ Notice that our definitions of $\mathscr{X}^{\bot}$ and $^{\bot}\mathscr{X}$
are not the same as \cite{BR} and \cite{BK}. Nevertheless, the following result indicates that our definitions of $\mathcal{GP}^{\bot}$ and $^{\bot}\mathcal{GI}$ coincide with those in \cite{BR} and \cite{BK}.

\begin{lemma}\label{21} The following statements hold.
\begin{enumerate}
\item[(i)] If $M\in \mathcal{GP}^{\bot},$ then $\Ext^{i}_{R}(Q,M)=0$ for all $i\geqslant 1$ and all
$Q\in \mathcal{GP}.$
\item[(ii)] If $M\in ^{\bot}\mathcal{GI},$ then $\Ext^{i}_{R}(M,E)=0$ for all $i\geqslant 1$ and all
$E\in \mathcal{GI}.$
\end{enumerate}
\end{lemma}

\begin{prf} (i) Let $Q$ be a Gorenstein projective $R$-module and $i\geqslant 2$. Then by the definition, there
exists an exact sequence $$\xymatrix{0\ar[r]&\widetilde{Q}\ar[r]&P_{i-2}\ar[r]&\ldots\ar[r]&P_{0}\ar[r]&Q\ar[r]&0}$$
of $R$-modules, where $P_{0},\ldots ,P_{i-2}$ are projective and $\widetilde{Q}$ is Gorenstein projective. Now for
any $M\in \mathcal{GP}^{\bot},$ as  $\widetilde{Q}$ is Gorenstein projective, one has $\Ext^{i}_{R}(Q,M)\cong \Ext^{1}_{R}(\widetilde{Q},M)=0.$

(ii) The proof is similar to the proof of (i), and so we leave it to the reader.
\end{prf}

We need the following lemma in the proof of Proposition \ref{26}.

\begin{lemma}\label{23} The following statements hold.
\begin{enumerate}
\item[(i)] We have $\mathcal{GP}\subseteq \{M\in R-\Mod|M \  \text{is}\  \mathcal{GP}-\text{pure projective}\}
\subseteq ^{\bot}(\mathcal{GP}^{\bot}).$
\item[(ii)] Every $\mathcal{GP}$-pure projective $R$-module of finite Gorenstein projective dimension is
Gorenstein projective.
\end{enumerate}
\end{lemma}

\begin{prf} (i) The left containment is trivial by the definition. Now, let $P$ be a $\mathcal{GP}$-pure
projective $R$-module and $M\in \mathcal{GP}^{\bot}$. One could have an exact sequence $${\bf X}=
\xymatrix{0\ar[r]& M\ar[r]& E\ar[r]& C\ar[r]& 0,}$$ of $R$-modules in which $E$ is injective. Then for
any $Q\in\mathcal{GP},$ as $\Ext^{1}_{R}(Q,M)=0,$ we deduce that the sequence $$\xymatrix{0\ar[r]&\Hom_{R}(Q,M)\ar[r]&\Hom_{R}(Q,E)\ar[r]&\Hom_{R}(Q,C)\ar[r]&0}$$ is exact. So,
${\bf X}$ is $\mathcal{GP}$-pure exact. Then, as $P$ is $\mathcal{GP}$-pure projective, $\Hom_{R}(P,{\bf X})$
is exact, and so the exact sequence $$\xymatrix{0\ar[r]&\Hom_{R}(P,M)\ar[r]&\Hom_{R}(P,E)\ar[r]&\Hom_{R}(P,C)\ar[r]&\Ext^{1}_{R}(P,M)
\ar[r]&0}$$ implies that $\Ext^{1}_{R}(P,M)=0.$ Thus $P\in\ ^{\bot}(\mathcal{GP}^{\bot}),$ as required.

(ii) Let $M$ be a $\mathcal{GP}$-pure projective $R$-module of finite Gorenstein projective dimension and 
let $\Gpd_RM$ denote the Gorenstein projective dimension of $M$. Clearly, we may assume that $M$ is nonzero. 
It is immediate from the definition of Gorenstein projective modules that every projective $R$-module is 
contained in $\mathcal{GP}^{\bot}$, and so by (i), one has $\Ext_R^i(M,P)=0$ for all projective $R$-modules 
$P$ and all $i>0$. Thus, by \cite[Theorem 3.1]{CFH} one has $$\Gpd_RM=\sup \{i\in \mathbb{N}_0|\Ext_R^i(M,P)
\neq 0 \  \ \text{for some projective} \ \text{$R$-module}  \ P\}=0.$$ 
\end{prf}

\begin{corollary}\label{24} Assume that $R$ is a coherent ring of finite Krull dimension. Then every $\mathcal{GP}$-pure projective $R$-module is Gorenstein flat.
\end{corollary}

\begin{prf} Let $\mathcal{GF}$ denote the class of Gorenstein flat $R$-modules. By \cite[Theorem 2.11]{EJL}, we
have $^{\bot}(\mathcal{GF}^{\bot})=\mathcal{GF}$. Since $R$ is coherent and $\dim R<\infty$, by \cite[Proposition 3.7]{CFH} every Gorenstein projective $R$-module is Gorenstein flat. Thus, $^{\bot}(\mathcal{GP}^{\bot})\subseteq ^{\bot}(\mathcal{GF}^{\bot}),$ and so Lemma \ref{23} (i) completes the proof.
\end{prf}

Next, we recall the definitions of special precovers and preenvelopes.

\begin{definition}\label{25}  Let $M$ be an $R$-module and $\mathscr{X}$ be a class of $R$-modules.
\begin{enumerate}
\item[(i)] A $\mathscr{X}$-precover $\varphi: X\to M$ is called {\it special} if $\varphi$ is surjective
and $\Ker \varphi\in \mathscr{X}^{\bot}$. Also,  A $\mathscr{X}$-\emph{preenvelope} $\psi:M\to X$ is called
{\it special} if $\psi$ is injective and $\Coker \psi\in ^{\bot}\mathscr{X}$.
\item[(ii)] A pair $(\mathscr{X},\mathscr{Y})$ of $R$-modules is called a \emph{cotorsion theory}
if $\mathscr{X}^{\bot}=\mathscr{Y}$ and $^{\bot}\mathscr{Y}=\mathscr{X}$. A cotorsion theory $(\mathscr{X},\mathscr{Y})$ is said to be \emph{complete} if every $R$-module has a special $\mathscr{X}$-precover; or equivalently if every $R$-module has a special $\mathscr{Y}$-preenvelope.
\end{enumerate}
\end{definition}

A class $\mathscr{X}$ of $R$-modules is called \emph{precovering} (respectively \emph{preenvelopeing}) if
every $R$-module possesses a $\mathscr{X}$-precover (respectively $\mathscr{X}$-preenvelope).

\begin{proposition}\label{26} Assume that every $R$-module in $^{\bot}(\mathcal{GP}^{\bot})$ has a Gorenstein projective precover (this is the case if for instance the class $\mathcal{GP}$ is precovering). Then $\mathcal{GP}$ coincides with the class of $\mathcal{GP}$-pure projective $R$-modules.
\end{proposition}

\begin{prf} Obviously, any Gorenstein projective $R$-module is $\mathcal{GP}$-pure projective. Assume that $P$
is a $\mathcal{GP}$-pure projective $R$-module. Then, Lemma \ref{23} (i) yields that $P\in ^{\bot}(\mathcal{GP}^{\bot}).$ Thus, by the hypothesis $P$ admits a Gorenstein projective precover. Hence, there
is an exact sequence $$\xymatrix{{\bf X}=0\ar[r]& K\ar@{^{(}->}[r]^{i}& G\ar[r]^{\varphi}& P\ar[r]& 0,}$$ in which $G$ is Gorenstein projective and $\Hom_{R}(\widetilde{G},{\bf X})$ is exact for all Gorenstein projective $R$-modules $\widetilde{G}.$ Equivalently, ${\bf X}$ is $\mathcal{GP}$-pure exact. So, as $P$ is $\mathcal{GP}$-pure projective,
the sequence
$$\xymatrix{0\ar[r]&\Hom_{R}(P,K)\ar[r]^{i_{*}}&\Hom_{R}(P,G)\ar[r]^{\varphi_{*}}&\Hom_{R}(P,P)\ar[r]&0}$$ is exact. Consequently, ${\bf X}$ splits.  Since, by \cite[Theorem 2.5]{Ho1} the class $\mathcal{GP}$ is closed under direct summands, we deduce that $P$ is Gorenstein projective.
\end{prf}

We end this section by the following useful remark.

\begin{remark}\label{36}
\begin{enumerate}
\item[(i)] Each of the theories of classical homology, pure homology, RD-pure homology and cyclically pure
homology is a relative pure homology that its corresponding class $\mathscr{X}$ is a class of finitely presented $R$-modules; see \cite{W}.
\item[(ii)] There is no class  $\mathscr{X}$ of finitely presented $R$-modules such that Gorenstein
homology is a $\mathscr{X}$-pure homology. To this end, let $R$ be an Artinian Gorenstein local ring which is not a principal ideal ring. (The ring $\mathbb{R}[x,y]/\langle x^2,y^2 \rangle$ is an instance of such a ring.) Assume in the contrary that there is a class $\mathscr{X}$ of finitely presented $R$-modules such that Gorenstein homology is a $\mathscr{X}$-pure homology. Then, every $R$-module is $\mathscr{X}$-pure projective. Note that every $R$-module is Gorenstein projective, because $R$ is assumed to be an Artinian Gorenstein local ring. As any pure exact sequence is also $\mathscr{X}$-pure exact, it follows that every $R$-module is pure projective. Hence, the global pure projective dimension of $R$ is zero. Now, \cite[Theorem 4.3]{Gr} yields that $R$ is a principal ideal ring. This contradicts our assumption on $R$.
\item[(iii)] Enochs, Jenda and L\'{o}pez-Ramos \cite[Theorem 2.12]{EJL} have proved that over a coherent ring
$R$, every $R$-module possesses a Gorenstein flat cover. Assume that $R$ is an arbitrary ring admitting a class $\mathscr{X}$ of $R$-modules such that the class of Gorenstein flat $R$-modules coincides with the class of $\mathscr{X}$-pure flat $R$-modules. Then \cite[Corollary 2.3]{ZD} implies that every $R$-module possesses a Gorenstein flat cover.
\item[(iv)] Assume that $R$ is such that $\mathcal{GP}$, the class of Gorenstein projective $R$-modules, is
precovering. Then for any given $R$-module $M$, we have a complex ${\bf Q}_{\bullet}=\cdots \to Q_{i}\overset{d_i}\to \cdots\overset{d_1}\to Q_{0}\to 0$ of Gorenstein projective $R$-modules and an $R$-homomorphism $\varphi:Q_0\to M$ such that the complex $\cdots \to Q_{i}\overset{d_i}\to \cdots\overset{d_1}\to Q_{0}\overset{\varphi}\to M\to 0$ is $\mathcal{GP}$-pure exact; see Definition \ref{35} below. We define $\Ext^{i}_{\mathcal{GP}}(M,N):= H^{i}(\Hom_R({\bf Q}_{\bullet},N))$ for all $R$-modules $N$ and all $i\geq 0$. Definition of $\Ext^{i}_{\mathcal{GP}}(M,N)$ is independent of the choose of ${\bf Q}_{\bullet}$. J{\o}rgensen
\cite[Corollary 2.13]{J} proved that over any  Noetherian ring with dualizing complex, the class $\mathcal{GP}$ is precovering. Then, Murfet and Salarian \cite[Theorem A.1]{MS} extended his result to Noetherian rings of finite Krull dimension. More recently, it is shown that if $R$ is coherent and every flat $R$-module has finite projective dimension, then $(\mathcal{GP},\mathcal{GP}^{\bot})$ is a complete cotorsion theory. This is part of the ongoing work \cite{BGH} which is also reported in \cite{Gi}.
\item[(v)] Let $R$ be such that $\mathcal{GI}$, the class of Gorenstein injective $R$-modules, is
preenvelopeing. Then for any given $R$-module $N$,  there is a complex ${\bf E}^{\bullet}=0\to E^{0}\overset{d^0}\to \cdots \to E^{i}\overset{d^i}\to \cdots$ of Gorenstein injective $R$-modules and an $R$-homomorphism $\psi:N\to E^0$ such that the complex  $0\to N\overset{\psi}\to E^{0}\overset{d^0}\to \cdots \to E^{i}\overset{d^i}\to \cdots$ is $\mathcal{GI}$-copure exact; see Definition \ref{35} below.  We define $\Ext^{i}_{\mathcal{GI}}(M,N):= H^{i}(\Hom_R(M,{\bf E}^{\bullet}))$ for all $R$-modules $M$ and all $i\geq 0$. This definition is well-defined.
Let $R$ be a Noetherian ring. Then \cite[Theorem 7.12]{K} yields that $(^{\bot}\mathcal{GI},\mathcal{GI})$
is a complete cotorsion theory. In particular,  $\mathcal{GI}$ is a preenvelopeing class.
\end{enumerate}
\end{remark}

\section{Gorenstein And Relative Pure Injectivity And Flatness}

In what follows, we denote the Pontryagin duality functor $\Hom_{\mathbb{Z}}(-,\mathbb{Q}/\mathbb{Z})$
by $(-)^{+}$. We summarize \cite[Lemma 2.2]{ZD} and \cite[Theorem 3.6]{Ho1} in the following lemma.

\begin{lemma}\label{31} Let $R$ be a coherent ring, $M$ an $R$-module and $\mathscr{X}$ a class of $R$-modules.
\begin{enumerate}
\item[(i)] $M$ is $\mathscr{X}$-pure flat if and only if $M^{+}$ is $\mathscr{X}$-pure injective.
\item[(ii)] $M$ is Gorenstein flat if and only if $M^{+}$ is Gorenstein injective.
\end{enumerate}
\end{lemma}

We record the following immediate corollary. 

\begin{corollary}\label{32} Let $R$ be a coherent ring. Assume that the class of $\mathcal{GP}$-pure injective $R$-modules coincides with the class of Gorenstein injective $R$-modules. Then the class of $\mathcal{GP}$-pure flat $R$-modules coincides with the class of Gorenstein flat $R$-modules.
\end{corollary}

The next result implies the Gorenstein injective and Gorenstein flat counterparts of Lemma \ref{23} (ii).

\begin{lemma}\label{33} The following statements hold.
\begin{enumerate}
\item[(i)] Every $\mathcal{GP}$-pure injective $R$-module of finite Gorenstein injective dimension is
Gorenstein injective.
\item[(ii)] Assume that $R$ is coherent. Every $\mathcal{GP}$-pure flat $R$-module of finite Gorenstein
flat dimension is Gorenstein flat.
\end{enumerate}
\end{lemma}

\begin{prf} (i) Let $M$ be a $\mathcal{GP}$-pure injective $R$-module of finite Gorenstein injective dimension.
In view of \cite[Theorem 2.15]{Ho1} and \cite[Lemma 3.5]{Ho2}, there exists a $\mathcal{GP}$-pure exact sequence $$\xymatrix{0\ar[r]& M\ar[r]^{f}& E\ar[r]& C\ar[r]&0,} \   \   (\dag)$$ where $E$ is Gorenstein injective.
Since $M$ is $\mathcal{GP}$-pure injective, we obtain the following exact sequence
$$\xymatrix{0\ar[r]& \Hom_{R}(C,M)\ar[r]&\Hom_{R}(E,M)\ar[r]^{f^{*}}&\Hom_{R}(M,M)\ar[r]&0.}$$
Hence $(\dag)$ splits, and so by \cite[Theorem 2.6]{Ho1}, $M$ is Gorenstein injective.

(ii) Let $M$ be a $\mathcal{GP}$-pure flat $R$-module of finite Gorenstein flat dimension. Since $M$ has a finite Gorenstein flat resolution, Lemma \ref{31} yields that $M^{+}$ is $\mathcal{GP}$-pure injective and it has finite Gorenstein injective dimension. Thus by (i), we obtain that $M^{+}$ is Gorenstein injective. By using Lemma \ref{31} again, we deduce that $M$ is Gorenstein flat.
\end{prf}

Let $\mathscr{X}$ be a class of $R$-modules and $\mathbf{Y}:=0\to L\to M\to N\to 0$ an exact sequence of
$R$-modules. We say  $\mathbf{Y}$ is $\mathscr{X}$-{\it copure exact} if  $\Hom_{R}({\bf Y},V)$ is exact for all
$V\in  \mathscr{X}.$ In what follows, we will use the following standard fact:

\begin{lemma}\label{34} Let $\mathscr{X}$ be a class of $R$-modules and $\mathbf{X}=\cdots\to X_{i+1}\to\hspace{-6mm}^{\vspace{-1mm}^{d_{i+1}}} \hspace{1mm}X_{i} \to\hspace{-4mm}^{\vspace{-1mm}^{d_{i}}}\hspace{1mm} X_{i-1}\to\cdots$ an exact complex of $R$-modules.
For each $i\in \BZ$, set $\mathbf{X_{i}}:=0\to \im d_{i+1}\hookrightarrow\hspace{-3mm}^{\vspace{-1mm}^{}}
\hspace{1mm} X_{i}\to\hspace{-4mm}^{\vspace{-1mm}^{}}\hspace{2mm} \im d_{i}\to 0.$ Then
\begin{enumerate}
\item[(i)] $\Hom_{R}(U,\mathbf{X})$ is exact for all $U\in  \mathscr{X}$ if and only if $\mathbf{X_{i}}$ is
$\mathscr{X}$-pure exact for all $i\in \BZ.$
\item[(ii)] $\Hom_{R}({\bf X},V)$ is exact for all $V\in  \mathscr{X}$ if and only if $\mathbf{X_{i}}$ is
$\mathscr{X}$-copure exact for all $i\in \BZ$.
\end{enumerate}
\end{lemma}

\begin{definition}\label{35} Let $\mathscr{X}$ be a class of $R$-modules. An exact complex $\mathbf{X}$
of $R$-modules is said to be $\mathscr{X}$-{\it pure exact} (respectively $\mathscr{X}$-{\it copure exact})
if it satisfies the equivalent conditions of Lemma \ref{34} (i) (respectively Lemma \ref{34} (ii)).
\end{definition}

The following lemma will be used several times in the rest of the paper. Its second part is also followed by
\cite[Theorem 7.12 (4)]{K} when $R$ is assumed to be Noetherian.

\begin{lemma}\label{37} Let $\mathcal{P}$ and $\mathcal{I}$ be the classes of projective and injective $R$-modules, respectively. Then the following assertions hold.
\begin{enumerate}
\item[(i)] $\mathcal{GP}\cap \mathcal{GP}^{\bot}=\mathcal{P}.$
\item[(ii)] $\mathcal{GI}\cap ^{\bot}\mathcal{GI}=\mathcal{I}.$
\end{enumerate}
\end{lemma}

\begin{prf} The proofs of (i) and (ii) are similar, and so we only prove (i).

Clearly, any projective $R$-module is Gorenstein projective. On the other hand, from the definition of
Gorenstein projective modules, it is immediate that $\Ext^i_R(Q,P)=0$ for all $Q\in \mathcal{GP}$, $P
\in \mathcal{P}$ and $i>0$. Hence, $\mathcal{P}\subseteq \mathcal{GP}\cap \mathcal{GP}^{\bot}$. Next,
let $M\in \mathcal{GP}\cap \mathcal{GP}^{\bot}.$ As $M$ is Gorenstein projective, there exists an exact
sequence $$\xymatrix{0\ar[r]& M\ar[r]& P\ar[r]& \widetilde{M}\ar[r]& 0} \  \  (*)$$ of $R$-modules, in
which $P$ is projective and $\widetilde{M}$ is Gorenstein projective. Since $M\in \mathcal{GP}^{\bot}$,
we have $\Ext_R^1(\widetilde{M},M)=0$, and so $(*)$ splits. Thus $M$ is projective.
\end{prf}

\begin{lemma}\label{38} Assume that $R$ is a Noetherian ring of finite Krull dimension. Then $ ^{\bot}\mathcal{GI}
\subseteq  \mathcal{GP}^{\bot}.$
\end{lemma}

\begin{prf} Let $M\in ^{\bot}\mathcal{GI}.$  By Remark \ref{36} (v), there exists a $\mathcal{GI}$-copure exact sequence $0\to M\overset{f}\to I\overset{g}\to C\to 0,$ where $I$ is Gorenstein injective and $C\in ^{\bot}\mathcal{GI}$. Since both $M$ and $C$ belong to $^{\bot}\mathcal{GI},$ it follows that $I\in ^{\bot}\mathcal{GI}.$ So, $I$ is injective by Lemma \ref{37} (ii). Next, as by Remark \ref{36} (iv) the pair $(\mathcal{GP},\mathcal{GP}^{\bot})$ is a complete cotorsion theory, we can find an exact sequence $0\to M\overset{h}\to N\overset{k}\to G\to 0,$ where $N\in \mathcal{GP}^{\bot}$ and $G\in \mathcal{GP}.$

Consider the following pushout diagram
$$\xymatrix{&0\ar[d]&0\ar[d]& & \\
 0\ar[r]&M\ar[d]_{h}\ar[r]^{f}&I\ar[d]_{\beta}\ar[r]^{g}&C\ar@{=}[d]\ar[r]&0\\
             0\ar[r]&N\ar[d]_{k}\ar[r]^{\alpha}&X\ar[d]\ar[r]^{\gamma}&C\ar[r]&0\\
                    &G\ar[d]\ar@{=}[r]&G\ar[d]& & \\
                    & 0 &0 & &}$$
where $X:=(N\oplus I)/S$ with $S=\{(h(m),-f(m))\ |\ m\in M \}$ and $\gamma((x,y)+S)=g(y)$ for all $(x,y)+S\in X.$
We intend to prove that the upper short exact sequence is $\mathcal{GP}$-pure exact. To this end, we have to show
that for every Gorenstein projective $R$-module $P$, the induced $R$-homomorphism $\Hom_R(P,I)\lo \Hom_R(P,C)$ is
surjective. Let $P$ be a Gorenstein projective $R$-module and $\varphi:P\lo C$ be an $R$-homomorphism. As $N\in \mathcal{GP}^{\bot},$ the lower short exact sequence is $\mathcal{GP}$-pure exact. So, there exists an $R$-homomorphism $\psi:P\to X$ such that $\gamma \psi=\varphi.$ Thus, we have the following commutative
diagram:
$$\xymatrix{
 0\ar[r]&M\ar[d]_{h}\ar[r]^{f}&I\ar[d]_{\beta}\ar[r]^{g}&C\ar@{=}[d]\ar[r]&0\\
             0\ar[r]&N\ar[r]^{\alpha}&X\ar[r]^{\gamma}&C\ar[r]&0\\
                    &&& P\ar@{-->}[ul]_{\psi}\ar[u]_{\varphi}&
                   & &}$$
Because $\beta:I\lo X$ is one-to-one and the $R$-module $I$ is injective, there is an $R$-homomorphism
$\rho:X\lo I$ such that $\rho \beta=\id_I$. Set $\theta:=\rho \psi$. Then $\theta \in \Hom_R(P,I),$ and it
can be easily verified that $g\theta=\varphi.$ Thus the upper short exact sequence is $\mathcal{GP}$-pure exact,
as required.

Now, let $Q$ be a Gorenstein projective $R$-module. As $I$ is injective, one has $\Ext^{1}_{R}(Q,I)=0.$ Hence,
we have the following exact sequence
$$\xymatrix{0\ar[r]&\Hom_{R}(Q,M)\ar[r]&\Hom_{R}(Q,I)\ar[r]&\Hom_{R}(Q,C)\ar[r]&\Ext^{1}_{R}(Q,M)\ar[r]&0.}$$
This implies that $\Ext^{1}_{R}(Q,M)=0,$ and so $M\in\ \mathcal{GP}^{\bot}.$
\end{prf}

The following corollary will be needed in the proof of Theorem \ref{41}.

\begin{corollary}\label{39} Assume that $R$ is a Noetherian ring of finite Krull dimension. Then every $R$-module
$M$ admits a $\mathcal{GI}$-copure exact complex $0\to M\overset{\psi}\to E^{0}\overset{d^0}\to \cdots \to E^{i}\overset{d^i}\to \cdots$ which is $\mathcal{GP}$-pure exact and each $E^i$ is Gorenstein injective.
\end{corollary}

\begin{prf}  Let $M$ be an $R$-module. By Remark \ref{36} (v), $M$ possesses a Gorenstein injective preenvelope $\psi:M\to E$ with $C:=\Coker \psi\in ^{\bot}\mathcal{GI}.$ Then, the sequence $$\xymatrix{{\bf X}=0\ar[r]& M\ar[r]^{\psi}&E\ar[r]^{\pi}& C\ar[r]& 0}$$ is $\mathcal{GI}$-copure exact. Let $Q$ be a Gorenstein projective $R$-module. We show that $\Hom_{R}(Q,{\bf X})$ is exact. Equivalently, we prove that $\pi_{*}:\Hom_{R}(Q,E)\to \Hom_{R}(Q,C)$ is surjective. Let $\alpha:Q\to C$ be an $R$-homomorphism. There exists an exact sequence $$\xymatrix{0\ar[r]&Q\ar[r]^{f}&P\ar[r]&\widetilde{Q}\ar[r]&0,}$$ where $P$ is projective and $\widetilde{Q}$ is Gorenstein projective.  Lemma \ref{38} implies that $C\in \mathcal{GP}^{\bot},$ and so $\Ext^{1}_{R}(\widetilde{Q},C)=0.$ Hence, $f^{*}:\Hom_{R}(P,C)\to \Hom_{R}(Q,C)$ is surjective, and so there
exists an $R$-homomorphism $\beta:P\to C$ such that $\alpha=f^{*}(\beta)=\beta f.$ Since $P$ is projective, we
have an $R$-homomorphism $g:P\to E$ making the following diagram commutative
$$\xymatrix{
& &P\ar[d]^{\beta}\ar@{-->}[dl]_{g}\\
&E \ar[r]_{\pi}&C
 }$$
Now, for the $R$-homomorphism $gf:Q\to E$, one has $$\pi_{*}(gf)=\pi(gf)=\beta f=\alpha.$$ By continuing the above argument and applying Lemma \ref{34}, we obtain a $\mathcal{GI}$-copure exact complex $0\to M\overset{\psi}\to E^{0}\overset {d^0}\to \cdots \to E^{i}\overset{d^i}\to \cdots$ which is $\mathcal{GP}$-pure exact and each $E^i$
is Gorenstein injective.
\end{prf}

The notion of virtually Gorenstein Artin algebras was introduced by Beligiannis and Reiten in \cite{BR}; see also
\cite{BK}. Next, we extend this notion to commutative rings.

\begin{definition}\label{40}  A Noetherian ring $R$ of finite Krull dimension is called \emph {virtually Gorenstein} if $\mathcal{GP}^{\bot}=^{\bot}\mathcal{GI}.$
\end{definition}

Next, we present a characterization of virtually Gorenstein rings. Since in its statement the phrase ``right
balanced'' is appeared, we recall  the meaning of this phrase here. Let $\mathcal{F}$ and $\mathcal{G}$ be two
classes of $R$-modules.  The functor $\Hom_R(-,\sim)$ is said to be \emph{right balanced} by $\mathcal{F}\times
\mathcal{G}$ if for every $R$-module $M$, there are complexes  $${\bf F_{\bullet}}=\cdots \to F_n\to F_{n-1}\to
\cdots \to F_0\to M\to 0$$ and
$${\bf G^{\bullet}}=0\to M\to G^0\to \cdots \to G^n\to G^{n+1}\to \cdots$$ in which $F_n\in \mathcal{F},
G^n\in \mathcal{G}$ for all $n\geq 0,$ such that for any $F\in \mathcal{F}$ and any $G\in \mathcal{G}$,
the two complexes $\Hom_R({\bf F_{\bullet}},G)$ and $\Hom_R(F,{\bf G^{\bullet}})$ are exact.

\begin{theorem}\label{41} Let $R$ be a Noetherian ring of finite Krull dimension and let $\mathcal{GP}$ and $\mathcal{GI}$ be the classes of Gorenstein projective and Gorenstein injective $R$-modules, respectively. The following are equivalent:
\begin{enumerate}
\item[(i)]  $\Hom_R(-,\sim)$ is right balanced by $\mathcal{GP}\times\mathcal{GI}$.
\item[(ii)] A short exact sequence ${\bf X}=0\to X_{1}\to X_{2}\to X_{3}\to 0$ of $R$-modules
is $\mathcal{GP}$-pure exact if and only if it is $\mathcal{GI}$-copure exact.
\item[(iii)] $R$ is virtually Gorenstein.
\end{enumerate}
\end{theorem}

\begin{prf} (i)$\Rightarrow$(ii) Assume that ${\bf X}$ is $\mathcal{GP}$-pure exact and $E$ is a Gorenstein
injective $R$-module. Then, by \cite[Theorem 8.2.3 (2)]{EJ2} we obtain the following exact sequence $$\xymatrix{0\ar[r]&\Hom_{R}(X_{3},E)\ar[r]&\Hom_{R}(X_{2},E)\ar[r]&\Hom_{R}(X_{1},E)\ar[r]&
\Ext^{1}_{\mathcal{GP}}(X_{3},E).}$$ As $E$ is Gorenstein injective, by \cite[Theorem 8.2.14]{EJ2} we conclude
that $$\Ext^{1}_{\mathcal{GP}}(X_{3},E)\cong \Ext^{1}_{\mathcal{GI}}(X_{3},E)=0.$$ Hence, ${\bf X}$ is
$\mathcal{GI}$-copure exact.  Similarly if ${\bf X}$ is $\mathcal{GI}$-copure exact, then by using
\cite[Theorem 8.2.5 (1)]{EJ2}, we can deduce that ${\bf X}$ is $\mathcal{GP}$-pure exact.

(ii)$\Rightarrow$(iii)  In view of Lemma \ref{38}, it is enough to show that $\mathcal{GP}^{\bot}\subseteq   ^{\bot}\mathcal{GI}.$ Let $M\in \mathcal{GP}^{\bot}.$ Since, by Remark \ref{36} (iv), $M$ admits a special
Gorenstein projective precover, we have an exact sequence $$\xymatrix{0\ar[r]& K\ar[r]& Q\ar[r]&M\ar[r]& 0,} \hspace{.5cm} (\ast)$$ in which $Q$ is Gorenstein projective and $K\in \mathcal{GP}^{\bot}$.  Then $(\ast)$ is  $\mathcal{GP}$-pure exact, and so by the assumption it is also $\mathcal{GI}$-copure. From $(\ast)$, we can see that $Q\in \mathcal{GP}^{\bot}.$ But, then Lemma \ref{37} (i) yields that $Q$ is projective. Let $E$ be a Gorenstein injective $R$-module. Since $\Ext^{1}_{R}(Q,E)=0,$  we deduce the
following exact sequence
$$\xymatrix{0\ar[r]&\Hom_{R}(M,E)\ar[r]&\Hom_{R}(Q,E)\ar[r]&\Hom_{R}(K,E)\ar[r]&\Ext^{1}_{R}(M,E)\ar[r]&0.}$$
But, the functor $\Hom_{R}(-,E)$ leaves $(\ast)$ exact, and so $\Ext^{1}_{R}(M,E)=0.$ Thus, $M\in ^{\bot}\mathcal{GI}.$

(iii)$\Rightarrow$(i) Let $M$ be an $R$-module. By an argument dual to the proof of Corollary \ref{39}, we
can construct a $\mathcal{GP}$-pure exact complex $\cdots \to Q_{i}\overset{d_i}\to \cdots\overset{d_1}\to Q_{0}\overset{\varphi}\to M\to 0$ which is $\mathcal{GI}$-copure exact and each $Q_i$ is Gorenstein projective.
Thus, in view of Corollary \ref{39}, it turns out that $\Hom_R(-,\sim)$ is right balanced by $\mathcal{GP}\times\mathcal{GI}$.
\end{prf}

In view of \cite[Theorem 8.2.14]{EJ2}, we record the following immediate corollary.

\begin{corollary}\label{42} Let $R$ be a virtually Gorenstein ring. Then  $\Ext^{i}_{\mathcal{GP}}(M,N)\cong \Ext^{i}_{\mathcal{GI}}(M,N)$ for all $R$-modules $M$ and $N$ and all $i\geq 0$.
\end{corollary}

The next result provides an answer to the main question of this investigation in the case the ground ring is
Noetherian of finite Krull dimension.

\begin{theorem}\label{43} Let $R$ be a Noetherian ring of finite Krull dimension and $\mathcal{GP}$ denote the class
of Gorenstein projective $R$-modules. The following are equivalent:
\begin{enumerate}
\item[(i)] $R$ is virtually Gorenstein.
\item[(ii)] the classes of Gorenstein injective and $\mathcal{GP}$-pure injective $R$-modules are the same.
\item[(iii)] Gorenstein homology is a $\mathcal{GP}$-pure homology.
\end{enumerate}
\end{theorem}

\begin{prf} (i)$\Rightarrow$(ii) Let $E$ be a $\mathcal{GP}$-pure injective $R$-module. Then $E$ has a Gorenstein injective preenvelope $\psi:E\to \widetilde{E}$ with $C:=\Coker \psi\in ^{\bot}\mathcal{GI};$ see Remark \ref{36} (v). One can easily see that
the sequence $$\xymatrix{{\bf X}=0\ar[r]&E\ar[r]^{\psi}&\widetilde{E}\ar[r]&C\ar[r]&0 }$$ is $\mathcal{GI}$
-copure exact, and so by Theorem \ref{41} it is also $\mathcal{GP}$-pure exact. Hence, as $E$ is $\mathcal{GP}$
-pure injective, the sequence $\Hom_{R}({\bf X},E)$ is exact, and so ${\bf X}$ splits. Thus $E$ is a direct
summand of $\widetilde{E},$ and so it is Gorenstein injective by \cite[Theorem 2.6]{Ho1}.

Conversely, let $E$ be a Gorenstein injective $R$-module. Let ${\bf X}=0\to X_{1}\to X_{2}\to X_{3}\to 0$ be
a $\mathcal{GP}$-pure exact sequence. Our assumptions on $R$ implies that the class $\mathcal{GP}$ is precovering,
and so by \cite[Theorem 8.2.3 (2)]{EJ2} we have the following exact sequence
$$\xymatrix{0\ar[r]&\Hom_{R}(X_{3},E)\ar[r]&\Hom_{R}(X_{2},E)\ar[r]&\Hom_{R}(X_{1},E)\ar[r]&
\Ext^{1}_{\mathcal{GP}}(X_{3},E)}$$ Now, Corollary \ref{42} yields that $$\Ext^{1}_{\mathcal{GP}}(X_{3},E)\cong\Ext^{1}_{\mathcal{GI}}(X_{3},E)=0,$$
and so $\Hom_{R}({\bf X},E)$ is exact. Thus, $E$ is $\mathcal{GP}$-pure injective.

(ii)$\Rightarrow$(i) By Lemma \ref{38}, it suffices to show that $\mathcal{GP}^{\bot}\subseteq ^{\bot}\mathcal{GI}.$ Let $M\in \mathcal{GP}^{\bot}$. Then, as by Remark \ref{36} (iv) $M$ has a special Gorenstein projective precover,
we have a $\mathcal{GP}$-pure exact sequence $$\xymatrix{0\ar[r]& K\ar[r]& P\ar[r]& M\ar[r]& 0,} \  \ (*)$$
where $P$ is Gorenstein projective and $K\in \mathcal{GP}^{\bot}$. As $M$ and $K$ belong to $\mathcal{GP}^{\bot}$, we deduce
that $P\in \mathcal{GP} \cap \mathcal{GP}^{\bot}$. Hence, $P$ is projective by Lemma \ref{37} (i). Since
every Gorenstein injective $R$-module is $\mathcal{GP}$-pure injective, it follows that $(*)$ is also $\mathcal{GI}$-copure exact. Let $E$ be a Gorenstein injective $R$-module. Then the functor $\Hom_{R}(-,E)$
leaves $(*)$ exact and $\Ext^{1}_{R}(P,E)=0.$ Thus, from the exact sequence
$$\xymatrix{0\ar[r]&\Hom_{R}(M,E)\ar[r]&\Hom_{R}(P,E)\ar[r]&\Hom_{R}(K,E)\ar[r]
&\Ext^{1}_{R}(M,E)\ar[r]&0,}$$ we conclude that $\Ext^{1}_{R}(M,E)=0.$ Equivalently, $M\in\ ^{\bot}\mathcal{GI}.$

(ii)$\Longleftrightarrow$(iii) follows in view of Remark \ref{36} (iv), Proposition \ref{26}, and Corollary \ref{32}.
\end{prf}

We conclude the paper with the following example which, in particular shows that there exists a ring such that the functor $\Hom_R(-,\sim)$ is not right balanced by $\mathcal{GP}\times\mathcal{GI}$.

\begin{example}\label{44}
\begin{enumerate}
\item[(i)] By \cite[Remarks 11.2.3 and 11.5.10]{EJ2}, every finite Krull dimensional Gorenstein ring is
virtually Gorenstein.
\item[(ii)] Let $(R,\fm)$ be any local ring with $\fm^2=0$. Then by \cite[Proposition 6.1 and Remark 6.5]{IK}
the only Gorenstein projectives are the projectives and the only Gorenstein injectives are the injectives.
Thus both classes $\mathcal{GP}^{\bot}$ and $^{\bot}\mathcal{GI}$ are equal to the class of all $R$-modules,
and so $R$ is virtually Gorenstein. So, we have plenty examples of virtually Gorenstein local rings which 
are not Gorenstein.
\item[(iii)] Let $k$ be a field. Then $R:=k[x,y,z]/\langle x^{2},yz,y^{2}-xz,z^{2}-yx\rangle$ is not
virtually Gorenstein by \cite[Proposition 4.3]{BK}. Hence, $\Hom_R(-,\sim)$ is not right balanced by $\mathcal{GP}\times\mathcal{GI}$.
\end{enumerate}
\end{example}

\subsection*{Acknowledgements} We thank the referee for her/his valuable comments. Also, we thank Srikanth
Iyengar for pointing out to us Example \ref{44} (ii).

\end{document}